\documentclass{amsart}
\usepackage{mathrsfs}
\usepackage{graphicx}
\usepackage{amsmath}
\usepackage{amssymb}
\usepackage{amsfonts}
\usepackage{amscd}
\usepackage{latexsym}
\usepackage{graphics}
\usepackage[all]{xy}

 \newtheorem{theorem}{Theorem}[section]
  \newtheorem{prop}[theorem]{Proposition}
  \newtheorem{cor}[theorem]{Corollary}
  \newtheorem{lemma}[theorem]{Lemma}
  \newtheorem{rmk}[theorem]{Remark}

  \newtheorem{thm}[theorem]{Theorem}

\def\det{\operatorname{det}}

\newcommand{\cal}{\mathcal}

\newcommand{\lra}{\longrightarrow}

\newcommand{\mt}{\mapsto}
\newcommand{\h}{\simeq}
\renewcommand{\sc}{\textsc}
\newcommand{\mb}{\mathbb}
\newcommand{\de}{\delta}
\renewcommand{\t}{\theta}

\newcommand{\f}{\frac}
\renewcommand{\l}{\left}
\renewcommand{\r}{\right}
\renewcommand{\lg}{\langle}
\newcommand{\rg}{\rangle}
\newcommand{\be}{\begin{equation}}
\newcommand{\ee}{\end{equation}}
\newcommand{\bce}{\begin{center}}
\newcommand{\ece}{\end{center}}
\newcommand{\bq}{\begin{eqnarray}}
\newcommand{\eq}{\end{eqnarray}}
\newcommand{\n}{\nonumber}
\renewcommand{\v}{\varepsilon}
\newcommand{\ra}{\rightarrow}

\renewcommand{\a}{\alpha}
\renewcommand{\b}{\beta}
\renewcommand{\rm}{\textrm}
\newcommand{\q}{\quad}

\newcommand{\la}{\lambda}

\newcommand{\g}{\gamma}
\newcommand{\G}{\Gamma}

\renewcommand{\i}{\infty}
\renewcommand{\c}{\cdot}
\newcommand{\cs}{\cdots}
\newcommand{\ba}{\begin{array}}
\newcommand{\ea}{\end{array}}
\renewcommand{\o}{\overline}
\newcommand{\om}{\omega}
\newcommand{\Om}{\Omega}

\renewcommand{\ss}{\subset}
\newcommand{\bpm}{\begin{pmatrix}}
\newcommand{\epm}{\end{pmatrix}}
\newcommand{\bbm}{\begin{bmatrix}}
\newcommand{\ebm}{\end{bmatrix}}
\newcommand{\ld}{\ldots}

\newcommand{\op}{\oplus}
\newcommand{\vp}{\varphi}
\newcommand{\ot}{\otimes}
\newcommand{\hr}{\hookrightarrow}

\newcommand{\bop}{\bigoplus}

\newcommand{\wtl}{\widetilde}

\newcommand{\bs}{\backslash}

\newcommand{\vt}{\vartheta}

\numberwithin{equation}{section}

\begin{document}

\title{Theta functions and arithmetic quotients of loop groups}

\author{Dongwen Liu}
\address{Department of Mathematics, Hong Kong University of Science and Technology}
\curraddr{Department of Mathematics, Hong Kong University of Science and Technology,
Clear Water Bay, Kowloon, Hong Kong}
\email{ldxab@ust.hk}


\date{\today}


\begin{abstract}
In this paper we observe that isomorphism classes of certain metrized
vector bundles over $\mb{P}^1_{\mb{Z}}-\{0,\i\}$ can be
parameterized by arithmetic quotients of loop groups. We construct
an asymptotic version of theta functions, which are defined on these
quotients. Then we prove the convergence and extend the theta functions to loop symplectic groups. We interpret them as sections of line bundles over an infinite dimensional torus,
discuss the relations with loop Heisenberg groups, and give an asymptotic multiplication formula.
\end{abstract}

\maketitle

\section{Introduction}

Let $E$ be a vector bundle of rank $n$ on the affine scheme $X:=\mb{P}^1_{\mb{Z}}-\{0,\i\}=\rm{Spec}R,$ where $R=\mb{Z}[t,t^{-1}].$
Then the global section $H^0(E)$ of $E$ is isomorphic to $R^n.$ Fix $q\in\mb{R},$ $q>1.$ If we choose a trivialization $H^0(E)\h R^n$, then we may
define an inner product on $H^0(E)$ via the following inner product $(,)_1$ on $R^n$:
\be\label{in}
(w(t),v(t))_1=\int_{q^{-1}S^1}w(t)^T \overline{v(t)} dt,
\ee
where $q^{-1}S^1=\{t\in\mb{C}: |t|=q^{-1}\},$ $\overline{v(t)}$ denotes the complex conjugate, and the supscript $^T$ is the transpose (here we write $w(t), v(t)$ as column vectors). The Lebesgue measure on $q^{-1}S^1$ is normalized such that meas$(q^{-1}S^1)=1.$ For any $g\in G:=GL_n(\mb{R}[t,t^{-1}]),$ we define
\be\label{ing}
(w,v)_g:=(g^{-1}w, g^{-1}v)_1.
\ee
We extend these inner products to $\mb{R}[t,t^{-1}]^n=R^n\ot_{\mb{Z}}\mb{R}$ by the same formulas.

Our starting point is to introduce the arithmetic quotient (see Section 2 for precise details)
\[
{\cal Q}_{n,q}=\wtl{\G}\bs\wtl{G}/\wtl{K},
\]
where $\G,$ $K$ are certain ``discrete'' and ``maximal compact" subgroup of $G$ respectively, $\wtl{G}$ is a central extension of $G$ by $\mb{R}^\times$, and $\wtl{\G}$, $\wtl{K}$ are double covers of $\G$, $K$. The key observation (Theorem \ref{moduli}) is that ${\cal Q}_{n,q}$ classifies above metrized vector bundles on $X$ with an additional datum called \textit{covolume theory}. Roughly speaking, a covolume theory, denoted by $c$, will assign a positive number $c(L)$ to each lattice $L$ of $H^0(E)$ in such a way that the information of relative covolumes between different lattices will be captured. Via this moduli interpretation of ${\cal Q}_{n,q}$, for any $\wtl{g}\in\wtl{G}$ we may associate a triple $(E,(,)_g,c)$. Then we construct an asymptotic version of theta functions on $\wtl{G}$
\[
\vt(\wtl{g}):=\lim_L c(L)\sum_{v\in L}e^{-\pi(v,v)_g},
\]
where $L$ runs over all lattices $L$ of $H^0(E)$, and the limit is taken with respect to inclusions of lattices. The interesting point of our construction is the (tautological!) fact that $\vt$ is defined on ${\cal Q}_{n,q}$. In other words, we obtain an automorphic function on $\wtl{G}$.

The most technical problem is the convergence of the limit defining $\vt(\wtl{g}).$ Using a variant of the lemmas in \cite{LZ} and \cite{Zhu}, together with some elementary Fourier analysis, in particular the Poisson summation formula, we establish in Section 3 the uniform convergence of $\vt(\wtl{g})$ for $\wtl{g}$ varying in certain Siegel subset of $\wtl{G}.$ In the course of the proof, the use of Iwasawa decompositions for loop groups simplify our considerations. Motivated by the work of Y. Zhu \cite{Zhu} on Weil representations, we also extend the theta functions to loop symplectic groups in Section 4. In all the above, the reparametrization $t\mt qt$ for the variable $t$, plays an important role for the convergence result. Note that it is also a crucial ingredient in the proof of the convergence of the Eisenstein series (see \cite{Gar}, \cite{L}) on loop groups.

We try to give some interpretations for our theta functions from representation-theoretic and geometric point of view. D. Mumford's idea \cite{Mum} is to bring together three ways of viewing theta functions:

(a) as holomorphic functions in the vector/period matrix,

(b) as matrix coefficients of a representation of the Heisenberg/metaplectic groups,

(c) as sections of line bundles on abelian varieties/moduli space of abelian varieties.

In Section 5, we interpret the theta functions as global sections of line bundles over an infinite dimensional torus, discuss the action of loop Heisenberg groups, and give an asymptotic multiplication formula, which is analogous to the theta relations given by Mumford in \cite{Mum}.

Inspirations from Kapranov's work \cite{Kap} are indispensable for the writing of this paper, especially for introducing the notion of covolume theory. We also use the notions of \textit{semiinfinite Grassmannian} and \textit{dimension theory} in \cite{Kap} to write our multiplication formula in Section 5 in a more conceptual way.

It should not be too surprising to get some generalizations of part of this paper to more general framework, say, to the situation of higher local fields and adelic spaces, or the so-called $C_n$-categories (see \cite{O}, \cite{OP}, \cite{OP2} and \cite{P2}). For instance, Parshin \cite{P} considered Heisenberg groups associated with objects from $C_0^{f.g.}$, the category of finitely generated abelian groups.  The loop Heisenberg groups introduced in Section 5 should be viewed as Heisenberg groups associated with objects from $C_1^{f.g.}$, the category of filtered abelian groups with finitely generated quotients.

{\bf Notations} For an abelian group $A$, we write $A_\mb{R}=A\ot_\mb{Z}\mb{R}.$ For a finite dimensional real vector space $V$, write $V^\ast=\rm{Hom}_\mb{R}(V,\mb{R}).$ Let
 $(,): V\times V^\ast\ra\mb{R}$ be the canonical pairing. If $L$ is a lattice of $V$, denote by $L^\ast\ss V^\ast$ the dual lattice with respect to the bi-character $e^{-2\pi i (,)}$. Let ${\cal S}(V)$ be the space of Schwartz functions on $V.$ If a Haar measure on $V$ has been chosen, then define the Fourier transform $\mathscr{F}: {\cal S}(V)\ra{\cal S}(V^\ast)$ by
\[
\mathscr{F}f(y)=\int_V f(x)e^{-2\pi i( x, y)}dx.
\]
The Haar measure on $V^\ast$ is determined by the condition that $\mathscr{F}$ is an isometry.

For a positive integer $l$, $\mb{Z}_l$ stands for $\mb{Z}/l\mb{Z}$ (in contrast to the usual meaning of $l$-adic integers for prime $l$), $U_l$ is the group of $l$-th roots of unity. The pairing $\mb{Z}_l\times U_l\ra U_l\ss S^1$ gives the Pontryagin duality. Denote $\zeta_l=e^{2\pi i/l}.$

For $d\in\mb{Z}$ we write $L^d_\mb{Z}=t^d\mb{Z}[t^{-1}]^n$, $L^d_\mb{R}=t^d\mb{R}[t^{-1}]^n=L^d_\mb{Z}\ot_\mb{Z}\mb{R}.$

For an invertible matrix $g$ write $g^\ast=(g^{-1})^T.$

\section{Moduli interpretations of arithmetic quotients}

Our purpose in this section is to give a moduli interpretation of arithmetic quotients of loop groups, which in turn classify the vector bundles
with metrics described in the introduction. We start from the loop group $GL_n(\mb{R}[t,t^{-1}]).$ It is well-known that this group has a central extension $\wtl{GL}_n(\mb{R}[t,t^{-1}])$ which we shall now define.

A \textit{lattice} $L$ of a
free $\mb{R}[t,t^{-1}]$-module $V$ of rank $n$ is a free
$\mb{R}[t^{-1}]$-submodule of rank $n$. In other words $L$ is an
$\mb{R}[t^{-1}]$-span of a basis of $V$. Any two lattices $L_1$, $L_2$ in
$V$ are \textit{commensurable}, in the sense that the quotients
$L_1/(L_1\cap L_2)$ and $L_2/(L_1\cap L_2)$ are finite-dimensional
over $\mb{R}$. For example any lattice in $\mb{R}[t,t^{-1}]^n$ is commensurable
with $\mb{R}[t^{-1}]^n.$

Let $L^0_\mb{R}$ be the lattice $\mb{R}[t^{-1}]^n,$ and $g\in
GL_n(\mb{R}[t,t^{-1}]).$ Since $gL^0_\mb{R}/(L^0_\mb{R}\cap
gL^0_\mb{R})$ and $L^0_\mb{R}/(L^0_\mb{R}\cap gL^0_\mb{R})$ are
finite-dimensional vector spaces over $\mb{R}$, we can talk about
their top wedge powers. Let $\det(L^0_\mb{R},gL^0_\mb{R})$ be the
tensor product
\[
\wedge^{top} (gL^0_\mb{R}/L^0_\mb{R}\cap gL^0_\mb{R})\ot_\mb{R} \wedge^{top} (L^0_\mb{R}/L^0_\mb{R}\cap
gL^0_\mb{R})^{-1},
\]
where $(\c)^{-1}=\rm{Hom}_\mb{R}(\c,\mb{R})$ denotes the dual vector space.
Then $\det(L^0_\mb{R},gL^0_\mb{R})$ is a one-dimensional vector space over $\mb{R}.$ Let $\det(L^0_\mb{R},gL^0_\mb{R})^\times$ be the set of nonzero vectors in $\det(L^0_\mb{R},gL^0_\mb{R})$, which form an $\mb{R}^\times$-torsor. Now define the
group
\be\label{lpr}
\wtl{GL}_n(\mb{R}[t,t^{-1}])=\{(g, \om_g)| g\in GL_n(\mb{R}[t,t^{-1}]),
\om_g\in \det(L^0_\mb{R}, gL^0_\mb{R})^\times\}.
\ee
The multiplication in the group is given by
\[
(g,\om_g)(h,\om_h)=(gh, \om_g\wedge g\om_h ),
\]
where $g\om_h$ is image of $\om_h$ under the natural map
\[
\det(L^0_\mb{R},
hL^0_\mb{R})\stackrel{g}{\lra} \det (gL^0_\mb{R}, ghL^0_\mb{R}),
\]
 and $\om_g\wedge g\om_h$
is defined via the isomorphism
\[
\det(L^0_\mb{R},gL^0_\mb{R})\wedge \det(gL^0_\mb{R},
ghL^0_\mb{R})\lra \det (L^0_\mb{R}, ghL^0_\mb{R}).\]
It turns out that $\wtl{GL}_n(\mb{R}[t,t^{-1}])$ is a central extension of $GL_n(\mb{R}[t,t^{-1}])$ by $\mb{R}^\times$, then we have the short exact sequence
\be\label{ext}
1\lra \mb{R}^\times\lra \wtl{GL}_n(\mb{R}[t,t^{-1}])\lra GL_n(\mb{R}[t,t^{-1}])\lra 1.
\ee

\

We are concerned with two subgroups of $GL_n(\mb{R}[t,t^{-1}])$ together with their extensions. The first one is $GL_n(\mb{Z}[t,t^{-1}])$, which plays the role of a discrete subgroup, or more precisely the so-called arithmetic subgroup. This group acts on $R^n=\mb{Z}[t,t^{-1}]^n.$ Similarly as above, let us define a lattice in $R^n$ to be a free $\mb{Z}[t^{-1}]$-submodule of rank $n$. Let $L^0_\mb{Z}=\mb{Z}[t^{-1}]^n$. For any $g\in GL_n(\mb{Z}[t,t^{-1}])$, $gL^0_\mb{Z}/(L^0_\mb{Z}\cap gL^0_\mb{Z})$ and $L^0_\mb{Z}/(L^0_\mb{Z}\cap gL^0_\mb{Z})$ are free $\mb{Z}$-modules of finite rank, and we let $\det(L^0_\mb{Z},gL^0_\mb{Z})$ be the rank one free $\mb{Z}$-module
\[
\wedge^{top} (gL^0_\mb{Z}/L^0_\mb{Z}\cap gL^0_\mb{Z})\ot_\mb{Z} \wedge^{top} (L^0_\mb{Z}/L^0_\mb{Z}\cap
gL^0_\mb{Z})^{-1}.
\]
Note that there is a natural isomorphism
\be \label{zr}
\vp:\q \det(L^0_\mb{Z},gL^0_\mb{Z})\ot_\mb{Z}\mb{R}\h \det(L^0_\mb{R},gL^0_\mb{R}).
\ee
Define the group
\be\label{lpz}
\wtl{GL}_n(\mb{Z}[t,t^{-1}])=\{(g,\om_g)| g\in GL_n(\mb{Z}[t,t^{-1}]),
\om_g\rm{ a basis of }\det(L^0_\mb{Z}, gL^0_\mb{Z})\}
\ee
which is a double cover of $GL_n(\mb{Z}[t,t^{-1}]),$ i.e. there is a short exact sequence
\[
1\lra \mb{Z}_2\lra \wtl{GL}_n(\mb{Z}[t,t^{-1}])\lra GL_n(\mb{Z}[t,t^{-1}])\lra 1.
\]
Via (\ref{zr}) we have an embedding $\wtl{GL}_n(\mb{Z}[t,t^{-1}])\hr \wtl{GL}_n(\mb{R}[t,t^{-1}])$, which is compatible with the
canonical embedding $\mb{Z}_2\hr \mb{R}^\times.$ In other words we have the commutative diagram
\begin{displaymath}
\xymatrix
{1 \ar[r] & \mb{Z}_2 \ar[d] \ar[r] & \wtl{GL}_n(\mb{Z}[t,t^{-1}]) \ar[d] \ar[r] & GL_n(\mb{Z}[t,t^{-1}]) \ar[r] \ar[d] & 1\\
1 \ar[r] & \mb{R}^\times  \ar[r] & \wtl{GL}_n(\mb{R}[t,t^{-1}]) \ar[r] & GL_n(\mb{R}[t,t^{-1}]) \ar[r]  & 1}
\end{displaymath}
This construction of double cover is quite natural (see \cite{Kap} $\S$2.0).

\

Recall (\ref{ing}) that for any $g\in GL_n(\mb{R}[t,t^{-1}])$ we have an inner product $(,)_g$ on $\mb{R}[t,t^{-1}]^n.$ Write
$\|\c\|_g$ for the induced norm. The following lemma is an easy exercise.

\begin{lemma}
Fix $g$ and for any subspace $H$ of $\mb{R}[t,t^{-1}]^n$ write
$\o{H}$ for the completion with respect to $\|\c\|_g.$ Then for any
two lattices $L_1$, $L_2$ of $\mb{R}[t,t^{-1}]^n,$
\[
\o{L}_1\cap\o{L}_2=\o{L_1\cap L_2},\q \f{L_i}{L_1\cap L_2}\h \f{\o{L}_i}{\o{L}_1\cap\o{L}_2}, i=1,2.
\]
\end{lemma}

\

Using this lemma, we may induce a metric on $\det(L^0_\mb{R}, gL^0_\mb{R})$, which we still denote by $\|\c\|_g.$
The subgroup
\be\label{cpt}
K=\{g\in GL_n(\mb{R}[t,t^{-1}])|g(q^{-1}t)g(q^{-1}t^{-1})^T=1\}
\ee
plays the role of a maximal compact subgroup. In fact $K$ is the isometry group of $(,)_1.$ Hence
$\|\c\|_{gk}=\|\c\|_g$ for any $g\in GL_n(\mb{R}[t,t^{-1}])$, $k\in K.$ We remark that by definition there is a family
of group homomorphisms
\be\label{ev}
ev(t):\q K\lra U(n)
\ee
for $t\in q^{-1}S^1$, via evaluation maps. For fixed $k\in K$, $ev(\c)(k): q^{-1}S^1\ra U(n)$ is smooth.
The double cover $\wtl{K}$ of $K$ is defined by
\be\label{lpk}
\wtl{K}=\{(g,\om_g)|g\in K, \om_g\in \det(L^0_\mb{R}, gL^0_\mb{R})^\times, \|\om_g\|_1=1\}.
\ee
Similarly, we also have the following commutative diagram
\begin{displaymath}
\xymatrix
{1 \ar[r] & \mb{Z}_2 \ar[d] \ar[r] & \wtl{K} \ar[d] \ar[r] & K \ar[r] \ar[d] & 1\\
1 \ar[r] & \mb{R}^\times  \ar[r] & \wtl{GL}_n(\mb{R}[t,t^{-1}]) \ar[r] & GL_n(\mb{R}[t,t^{-1}]) \ar[r]  & 1}
\end{displaymath}

\

We introduce the arithmetic quotient
\be\label{quo}
{\cal Q}_{n,q}:=\wtl{GL}_n(\mb{Z}[t,t^{-1}])\bs \wtl{GL}_n(\mb{R}[t,t^{-1}])/\wtl{K}.
\ee
The main result of this section is the following moduli interpretation of ${\cal Q}_{n,q}.$ Let $E$ be a vector bundle of rank $n$ on $X$, and $\|\c\|$ be a metric on $E$ such that after a trivialization of $E$ it is induced by (\ref{ing}) for some $g.$ A \textit{covolume theory} of $(E,\|\c\|)$ is a rule $c$ which assigns a positive real number $c(L)$ to each lattice $L$ of $H^0(E)$ such that
\[
c(L')=c(L)cov_{\|\c\|}\l(\f{L'}{L}, \f{L_\mb{R}'}{L_\mb{R}}\r)
\]
whenever $L\ss L'$ are two lattices of $H^0(E).$ Here $cov_{\|\c\|}$
denotes the covolume induced by the given metric $\|\c\|$, and
$L_\mb{R}=L\ot_\mb{Z}\mb{R}.$ Note that once we know $c(L)$ for one
lattice $L$, then we know $c(L')$ for any lattice $L'.$ In fact,
take a lattice $L''\ss L, L'$ and define \be\label{rcov}
c(\|\c\|,L,L')=\f{cov_{\|\c\|}\l(\f{L'}{L''},
\f{L_\mb{R}'}{L_\mb{R}''}\r)}{cov_{\|\c\|}\l(\f{L}{L''},
\f{L_\mb{R}}{L_\mb{R}''}\r)}. \ee Then (\ref{rcov}) does not depend
on the choice of $L''$ and satisfies
\be\label{ra}c(L')=c(L)c(\|\c\|,L,L'). \ee We remark that all
covolume theories attached to $(E,\|\c\|)$ form an
$\mb{R}^\times_+$-torsor.

Consider all the triples $(E,\|\c\|,c)$ as above. We say that two
such triples $(E_i,\|\c\|_i, c_i)$, $i=1,2$ are isomorphic if there
is an isomorphism $\eta: E_1\ra E_2$ of vector bundles such that
$\eta$ induces an isometry $H^0(E_1)\ra H^0(E_2)$ and
$c_2(\eta(L_1))=c_1(L_1)$ for any lattice $L_1$ of $H^0(E_1)$.

\begin{thm}\label{moduli}
${\cal Q}_{n,q}$ classifies the isomorphism classes of all such triples $(E,\|\c\|,c).$
\end{thm}
\textit{Proof.} Given $(g,\om_g)\in\wtl{GL}_n(\mb{R}[t,t^{-1}])$, we
associate the triple $(R^n, \|\c\|_g, c)$ such that the covolume
theory $c$ is determined by $c(L^0_\mb{Z})=\|\om_g\|^{-1}_g.$

(1) Let $(\g,\om_\g)\in\wtl{GL}_n(\mb{Z}[t,t^{-1}])^n.$ Then
associated to $(\g g, \om_\g\wedge \g\om_g)$ we have the triple
$(R^n,\|\c\|_{\g g},c')$ such that $c'$ is determined by
\[
c'(L^0_\mb{Z})=\|\om_\g\wedge \g\om_g\|_{\g g}^{-1}.
\]
Consider the automorphism of trivial vector bundles $\eta_\g: R^n\ra
R^n, v\mt \g v$. This induces an isometry because we have
\[
\|\g v\|_{\g g}=\|g^{-1}v\|_1=\|v\|_g.
\]
Moreover we have \bq c'(\g L^0_\mb{Z})&=&c'(L^0_\mb{Z})c(\|\c\|_{\g
g},L^0_\mb{Z},\g L^0_\mb{Z})\n\\
&=& \|\om_\g\wedge \g\om_g\|_{\g g}^{-1}c(\|\c\|_{\g
g},L^0_\mb{Z},\g L^0_\mb{Z})\n\\
&=&\|\om_g\|_g^{-1}\|\om_\g\|_{\g g}^{-1}c(\|\c\|_{\g
g},L^0_\mb{Z},\g L^0_\mb{Z})\n\\ &=&
\|\om_g\|_g^{-1}=c(L^0_\mb{Z}).\n
 \eq
This proves that $\eta_\g$ induces an isomorphism between the
triples $(R^n,\|\c\|_g,c)$ and $(R^n,\|\c\|_{\g g},c').$

(2) Let $(k,\om_k)\in \wtl{K}.$ Then associated to $(gk,\om_g\wedge
g\om_k)$ we have the triple $(R^n,\|\c\|_g, c'')$ such that $c''$ is
determined by
\[
c''(L^0_\mb{Z})=\|\om_g\wedge g\om_k\|_g^{-1}.
\]
In this case the identity map induces an isometry, and
\[
c''(L^0_\mb{Z})=\|\om_g\|_g^{-1}\|\om_k\|_1^{-1}=\|\om_g\|_g^{-1}=c(L^0_\mb{Z}).
\]
Hence the triples $(R^n,\|\c\|_g,c)$ and $(R^n,\|\c\|_g,c'')$ are
isomorphic via the identity map.

We have constructed a map from ${\cal Q}_{n,q}$ to the set of
isomorphism classes of such triples. The inverse map is also clear
and it is easy to prove this is a bijection. \hfill$\Box$

\section{Theta functions on the arithmetic quotients}

Suppose we have a triple $(E,\|\c\|,c)$ as in the previous section. Via isomorphisms we may always assume that $E=R^n$ and $\|\c\|=\|\c\|_g$ for some
$g\in GL(\mb{R}[t,t^{-1}])$. Let $\wtl{g}=(g,\om_g)\in\wtl{GL}_n(\mb{R}[t,t^{-1}])$ such that $\wtl{g}$ gives the triple $(R^n,\|\c\|_g,c)$ (see the proof of Theorem \ref{moduli}). We define
\be\label{theta}
\vartheta(\wtl{g})=\lim_L\vt_L(\wtl{g}):=\lim_L c(L)\sum_{v\in L}e^{-\pi(v,v)_g},
\ee
where $L$ runs over all lattices of $R^n\h H^0(E)$, and the limit is taken with respect to the inclusions of lattices. To make the definition more transparent, consider the filtration of lattices
\[
\cs\ss L^d_\mb{Z}\ss L^{d+1}_\mb{Z}\ss\cs
\]
where $L^d_\mb{Z}:=t^dL^0_\mb{Z}.$
Then $\bigcup\limits_{d\in\mb{Z}}L^d_\mb{Z}=R^n$ and
$\vt(\wtl{g})$ is the value (assume the existence) which satisfies
the condition that for any $\v>0$, there exists $d\in\mb{Z}$ such
that \be\label{cond} \l|\vt_L(\wtl{g})-\vt(\wtl{g})\r|<\v \ee whenever $L$ is a lattice
containing $L^d_\mb{Z}.$ Of course, in this description we may
replace $\{L^d_\mb{Z},d\in\mb{Z}\}$ by any exhausting
filtration of lattices. Clearly $\vt(\wtl{g})$ only depends on the
isomorphism class of the triple, and from Theorem \ref{moduli} it is
a tautological fact that $\vt$ descends to a function on ${\cal
Q}_{n,q}$, i.e. $\vt$ is left invariant under
$\wtl{GL}_n(\mb{Z}[t,t^{-1}])$ and right invariant under $\wtl{K}.$

\

In this section we shall prove uniform convergence properties of $\vt(\wtl{g})$, for $\wtl{g}$ varying in certain ``compact" subset of $\wtl{GL}_n(\mb{R}[t,t^{-1}]).$ Equivalently, we are considering behaviors of $\vt([\wtl{g}])$ when $[\wtl{g}]$ lies in certain subset of
the arithmetic quotient ${\cal Q}_{n,q}.$ To this end we need the following Iwasawa decomposition:

\begin{lemma}\label{iwa}
Let $B=\{g(t)\in GL_n(\mb{R}[t^{-1}])|g(\i)\rm{ is upper triangular}\}.$ The sequence $($\ref{ext}$)$ splits over $B$. If $\wtl{B}$ is the preimage of $B$ under the projection \[\wtl{GL}_n(\mb{R}[t,t^{-1}])\ra GL_n(\mb{R}[t,t^{-1}]),\] then
\[
\wtl{GL}_n(\mb{R}[t,t^{-1}])=\wtl{B}\wtl{K}=\mb{R}^\times B\wtl{K}.
\]
\end{lemma}

Proof of this lemma is by using standard theory of Tits systems. We also need a lemma which is similar to \cite{Zhu} Lemma 4.8 and \cite{LZ} Lemma 3.3.

\begin{lemma}\label{est}
Let $V_1,\ld, V_m$ be finite dimensional real vector spaces, $V=V_1\op\cs\op V_m.$ Suppose $f(x_1,\ld,x_m)\in {\cal S}(V)$ is a Schwartz function which takes values in $\mb{R}$ and its partial Fourier transforms $F_if\geq 0$ $(i=1,\ld,m)$, i.e.
\[
\int_{V_i}f(x_1,\ld,x_m)\lg x_i, y_i\rg dx_i\geq 0,
\]
where $y_i$ is in the Pontryagin dual of $V_i.$ Let $L_i\ss V_i$ be a lattice $($i.e. a free $\mb{Z}$-module which spans $V_i$ over $\mb{R})$ and $L=L_1\op\cs\op L_m$. If
$U: V\ra V$ is a blockwise unipotent linear operator, i.e.
\[
UV_i\ss V_1\op\cs\op V_i,\q U|_{V_i}\equiv 1 \mod V_1\op\cs\op V_{i-1},
\]
then
\[
\sum_{k\in L}f(Uk)\leq\sum_{k\in L}f(k).
\]
\end{lemma}

We give some notations. Let $N=\{g(t)\in
GL_n(\mb{R}[t^{-1}])|g(\i)=1\}$, and $A$ be the group of diagonal
matrices in $GL_n(\mb{R})$, then $B=AN.$ For each $d\in\mb{Z}$, let
$L^d_\mb{R}=t^dL^0_\mb{R}$, $V_d=t^d\mb{R}^n$ and $L_d=t^d\mb{Z}^n$, then \be\label{decom}
L^d_\mb{R}=\bop_{i\leq d}V_i,\q\q L^d_\mb{Z}=\bop_{i\leq d}L_i.
\ee
For $d_1<d_2\in\mb{Z}$, let
\be\label{fin}
V_{d_1,d_2}=\bop_{d_1\leq i\leq d_2}V_i,\q L_{d_1,d_2}=\bop_{d_1\leq i\leq d_2}L_i.
\ee
 Note that elements of $N$ act on $L^d_\mb{R}$ as blockwise
unipotent linear operators.

Let us define the notion of Siegel subsets of
$\wtl{GL}_n(\mb{R}[t,t^{-1}])$ which is suitable for our purpose.
$GL_n(\mb{R}[t,t^{-1}])$ has an ind-scheme structure. Indeed for
each $i\in\mb{N}$ let $X_i$ be the subset of
$GL_n(\mb{R}[t,t^{-1}])$ consisting of elements whose entries are
contained in $V_{-i,i}.$ Then
$GL_n(\mb{R}[t,t^{-1}])=\bigcup\limits_{i\in\mb{N}}X_i$, and one may
choose an embedding \be\label{aff} \vp_i: X_i\hr
\mb{A}^{N(i)}(\mb{R}) \ee of affine subscheme, where $N(i)\in\mb{N}$
depends on $i$. For example we may take $N(i)=n^2(2i+1)+2ni+1.$ We
call $C\ss GL_n(\mb{R}[t,t^{-1}])$ a compact subset, if $C$ is
contained in $X_i$ for some $i$, and $\vp_i(C)\ss\mb{R}^{N(i)}$ is
compact in the real topology. It is easy to check that this notion
is well-defined, i.e. does not depend on the choice of (\ref{aff}).
Similarly we can define compact subsets of $B.$
 Recall from Lemma
\ref{iwa} the Iwasawa decomposition. We define a \textit{Siegel subset} of
$\wtl{GL}_n(\mb{R}[t,t^{-1}])$ to be a subset of the form $C_Z
C_B\wtl{K}$ where $C_Z,$ $C_B$ are compact subsets of
$\mb{R}^\times$ and $B$ respectively.

Now we are ready to state and prove the main result of this section.

\begin{thm}
The limit $($\ref{theta}$)$ defining $\vt(\wtl{g})$ exists and
converges uniformly for $\wtl{g}$ varying in any Siegel subset of
$\wtl{GL}_n(\mb{R}[t,t^{-1}])$.
\end{thm}
\textit{Proof.} Let $u\in N$, $a\in A$, and suppose that $\wtl{g}$ gives the triple $(R^n,\|\c\|_{au},c).$
Our strategy is to first show that the sequence
\be\label{sq}
\vt_d(\wtl{g}):=c(L^d_\mb{Z})\sum_{v\in L^d_\mb{Z}}e^{-\pi(v,v)_{au}}
\ee
is bounded, and then compare $\vt_L(\wtl{g})$ and $\vt_d(\wtl{g})$ when $L$ is a lattice containing $L^d_\mb{Z}$. More precisely we shall prove that
\be\label{com}
\vt_L(\wtl{g})=(1+o_d(1))\vt_d(\wtl{g})
\ee
where $o_d(1)\ra 0$ as $d\ra\i$, and the rate of convergence depends on $d$ but not on $L$. We proceed in 3 steps.

\sc{Step 1}: Boundedness of subsequence. The trick is reducing to the diagonals. Let us first prove that for each $d\in\mb{Z}$,
\be\label{red}
\sum_{v\in L^d_\mb{Z}}e^{-\pi(v,v)_{au}}\leq\sum_{v\in L^d_\mb{Z}}e^{-\pi(v,v)_a}.
\ee
Lemma \ref{est} does not apply directly in this infinite dimensional case. However, consider
\be \label{psum}
\sum_{v\in L_{d_1,d}}e^{-\pi(v,v)_{au}}
\ee
for $d_1<d.$ Let $\pi_1,$ $\pi_2$ be the obvious projections
\[
\pi_1: L^d_\mb{R}\ra V_{d_1,d},\q \pi_2: L^d_\mb{R}\ra L^{d_1-1}_{\mb{R}}.
\]
Since $V_{d_1,d}$ and $L^{d_1-1}_{\mb{R}}$ are orthogonal with respect to $(,)_1$, we obtain
\bq
(v,v)_{au}&=&(u^{-1}a^{-1}v, u^{-1}a^{-1}v)_1\n\\
&\geq& (\pi_1 u^{-1}a^{-1}v, \pi_1 u^{-1}a^{-1}v)_1.\n
\eq
Apply Lemma \ref{est} for $V_i,$ $a^{-1}L_i$, $i=d_1,\ld, d$, $f(v)=e^{-\pi(v,v)_1}$ and the operator $U=\pi_1u^{-1}: V_{d_1,d}\ra V_{d_1,d}$, it follows that
\[
(\ref{psum})\leq \sum_{v\in a^{-1}L_{d_1,d}}e^{-\pi(Uv, Uv)_1}\leq \sum_{v\in a^{-1}L_{d_1,d}}e^{-\pi(v,v)_1}=\sum_{v\in L_{d_1,d}}e^{-\pi(v,v)_a}.
\]
Let $d_1\ra-\i$ we get (\ref{red}).

On the other hand, we have the equality of covolumes
\be\label{redcov}
cov_{\|\c\|_{au}}\l(\f{L^d_\mb{Z}}{L^{d_1}_\mb{Z}}, \f{L^d_\mb{R}}{L^{d_1}_\mb{R}}\r)=cov_{\|\c\|_{a}}\l(\f{L^d_\mb{Z}}{L^{d_1}_\mb{Z}}, \f{L^d_\mb{R}}{L^{d_1}_\mb{R}}\r).
\ee Indeed, take a basis $\{v_i\}$ of $L_{d_1,d}.$ For $v\in L^d_\mb{R}$ denote by $v^\perp$ the orthogonal projection of $v$ to $\l(L^{d_1}_\mb{R}\r)^\perp$ with respect to
the metric $\|\c\|_{au}$. It is easy to check that
\[
v^\perp=au\pi_1(u^{-1}a^{-1}v)=auUa^{-1}v.
\]
Then
\bq
&&\rm{LHS of }(\ref{redcov})\n\\
&=&\sqrt{\det\l((v_i^\perp, v_j^\perp)_{au}\r)}\n\\
&=&\sqrt{\det\l((Ua^{-1}v_i, Ua^{-1}v_j)_1\r)}\n\\
&=&\sqrt{\det\l((a^{-1}v_i,a^{-1}v_j)_1\r)}\n\\
&=&\rm{RHS of }(\ref{redcov}).\n
\eq

Write $a=\rm{diag}\{a_1,\ld,a_n\}.$ Apply (\ref{red}) and (\ref{redcov}), we obtain that for $d>0,$
\bq
\vt_d(\wtl{g})&\leq& c(L^0_\mb{Z})cov_{\|\c\|_{a}}\l(\f{L^d_\mb{Z}}{L^{0}_\mb{Z}}, \f{L^d_\mb{R}}{L^{0}_\mb{R}}\r)\sum_{v\in L^d_\mb{Z}}e^{-\pi(v,v)_a}\n\\
&=& c(L^0_\mb{Z})q^{-n\f{d(d+1)}{2}}|\det a|^{-d}\prod^n_{j=1}\prod^d_{i=-\i}\sum_{m\in\mb{Z}}\exp\l(-\pi a_j^{-2}q^{-2i}m^2\r).\n
\eq
By the Poisson summation formula,
\[
\sum_{m\in\mb{Z}}\exp\l(-\pi a_j^{-2}q^{-2i}m^2\r)=|a_j|q^i\sum_{m\in\mb{Z}}\exp\l(-\pi a_j^2 q^{2i}m^2\r).
\]
Hence
\be\label{bound}
\vt_d(\wtl{g})\leq c(L^0_\mb{Z})\prod^n_{j=1}\prod^\i_{i=0}\sum_{m\in\mb{Z}}\exp\l(-\pi a_j^{-2}q^{2i}m^2\r)\prod^d_{i=1}\sum_{m\in\mb{Z}}\exp\l(-\pi a_j^{2}q^{2i}m^2\r),
\ee
and the right-hand-side is convergent as $d\ra\i$.

\sc{Step 2}: Comparison with the subsequence. Let $L$ be a lattice containing $L^d_\mb{Z}$, and let $L'$ be a complement of $L^d_\mb{Z}$ in $L$ such that $L'\ss t^{d+1}\mb{Z}[t].$ Similarly as above, for $v'\in L'_\mb{R}$ write $v'^\perp$ for the orthogonal projection of $v'$ to $(L^d_\mb{R})^\perp.$ Write $v'_0=v'-v'^\perp\in L^d_\mb{R}.$ Then
\bq
\label{tl}\vt_L(\wtl{g})&=&c(L)\sum_{v\in L}e^{-\pi(v,v)_{au}}\\
&=&c(L)\sum_{v\in L^d_\mb{Z}}\sum_{v'\in L'}e^{-\pi(v+v',v+v')_{au}}\n\\
&=&c(L)\sum_{v'\in L'}e^{-\pi(v'^\perp,v'^\perp)_{au}}\sum_{v\in L^d_\mb{Z}}e^{-\pi(v+v'_0,v+v'_0)_{au}}   \n
\eq
It is clear that $L'_\mb{R}\perp L^{d-2\de(u)-1}_\mb{R}$, where $\de(u)$ is the highest power of
$t^{-1}$ appearing in the entries of $u$. As a consequence we have $v_0'\in V_{d-2\de(u),d}.$
For $v\in L^d_\mb{R}$ let $v_\perp$ be the orthogonal projection of $v$ to $V_{d-2\de(u),d}^\perp$, and let $v_0=v-v_\perp\in V_{d-2\de(u),d}.$ Then
\bq
&&\sum_{v\in L^d_\mb{Z}}e^{-\pi(v+v'_0,v+v'_0)_{au}}\label{pt2}\\
&=&\sum_{v\in t^{d-2\de(u)-1}L^0_\mb{R}}e^{-\pi(v_\perp,v_\perp)_{au}}\sum_{v_1\in L_{d-2\de(u),d}}e^{-\pi
(v_1+v_0+v_0',v_1+v_0+v_0')_{au}}. \n
\eq
Let $(,)^\ast_{au}$ be the induced inner product on $V_{d-2\de(u),d}^\ast$. Choose the Haar measure on $V_{d-2\de(u),d}$ such that the covolume of $L_{d-2\de(u),d}$ equals one (i.e. choose the ordinary Lebesgue measure on each $V_i=t^i\mb{R}\h\mb{R}$ via $rt^i\mt r$). Then for fixed $x_0\in V_{d-2\de(u),d}$, we have the Fourier transform
\[
\mathscr{F}\l(e^{-\pi(x+x_0,x+x_0)_{au}}\r)(y)=\f{e^{2\pi i (x_0, y)-\pi(y,y)^\ast_{au}}}{cov_{\|\c\|_{au}}(L_{d-2\de(u),d}, V_{d-2\de(u),d})}.
\]
By the Poisson summation formula,
\[
\sum_{v_1\in L_{d-2\de(u),d}}e^{-\pi
(v_1+x_0,v_1+x_0)_{au}}=\sum_{v^\ast_1\in L^\ast_{d-2\de(u),d}}\f{e^{2\pi i( x_0, v_1^\ast)-\pi(v_1^\ast,v_1^\ast)^\ast_{au}}}{cov_{\|\c\|_{au}}(L_{d-2\de(u),d}, V_{d-2\de(u),d})}.
\]
Previous reasoning for (\ref{red}) implies that
\bq
&&\l|\sum_{v^\ast_1\in L^\ast_{d-2\de(u),d}-\{0\}}e^{2\pi i(x_0, v_1^\ast)-\pi(v_1^\ast,v_1^\ast)^\ast_{au}}\r|\n\\
&\leq &\sum_{v^\ast_1\in L^\ast_{d-2\de(u),d}-\{0\}}e^{-\pi(v_1^\ast,v_1^\ast)^\ast_{au}}\n\\
&\leq & \sum_{v^\ast_1\in L^\ast_{d-2\de(u),d}-\{0\}}e^{-\pi(v_1^\ast,v_1^\ast)^\ast_{a}}\n\\
&=&\prod^d_{i=d-2\de(u)}\sum_{m\in\mb{Z}}\exp\l(-\pi a_j^{2}q^{2i}m^2\r)-1\n\\
&\leq & \prod^\i_{i=d-2\de(u)}\sum_{m\in\mb{Z}}\exp\l(-\pi a_j^{2}q^{2i}m^2\r)-1\n\\
&=&o_d(1).\n
\eq
Therefore we have proved that
\[
\sum_{v_1\in L_{d-2\de(u),d}}e^{-\pi
(v_1+x_0,v_1+x_0)_{au}}=\f{1+o_d(1)}{cov_{\|\c\|_{au}}(L_{d-2\de(u),d}, V_{d-2\de(u),d})}
\]
where the rate of convergence of $o_d(1)$ does not depend on $x_0.$ Hence from (\ref{pt2}) we obtain
\be\label{pt22}
\sum_{v\in L^d_\mb{Z}}e^{-\pi(v+v'_0,v+v'_0)_{au}}=(1+o_d(1))\sum_{v\in L^d_\mb{Z}}e^{-\pi(v,v)_{au}}
\ee
where again the magnitude of $o_d(1)$ does not depend on $v_0'$. Similar arguments apply for the  summation over $L'$ in (\ref{tl}), and yield
\be\label{pt1}
\sum_{v'\in L'}e^{-\pi(v'^\perp,v'^\perp)_{au}}
=\f{1+o_d(1)}{cov_{\|\c\|_{au}}\l(\f{L}{L^d_\mb{Z}},\f{L_\mb{R}}{L^d_\mb{R}}\r)}
\ee
where the magnitude of $o_d(1)$ does not depend on $L'$, hence not on $L$. Combine (\ref{tl}), (\ref{pt22}) and (\ref{pt1}), together with the identity
\[
c(L)=c(L^d_\mb{Z})cov_{\|\c\|_{au}}\l(\f{L}{L^d_\mb{Z}},\f{L_\mb{R}}{L^d_\mb{R}}\r),
\]
we obtain (\ref{com}).

\sc{Step 3}: Conclusions. After the first two steps we have proved the existence of $\vt(\wtl{g}).$ If we restrict to a Siegel subset, then by definition we have in the previous settings that $c(L^0_\mb{Z})$ and $au$ vary in some compact subsets of $\mb{R}^\times_+$ and $B$ respectively, and moreover $\de(u)$ is bounded. Consequently, the bound (\ref{bound}) and the quantities $o_d(1)$ appearing in Step 2 are all uniform. \hfill$\Box$

\section{Generalization to loop symplectic groups}

Recall that the Siegel upper half-space $\mathscr{H}_n$ for $Sp_{2n}(\mb{R})$ is the set of all $n\times n$ complex symmetric matrices with positive definite imaginary part. $Sp_{2n}(\mb{R})$ acts on $\mathscr
{H}_n$ via linear fractional transformations. $L^2(\mb{R}^n)$ is a model for the Weil representation of $Sp_{2n}(\mb{R})$, and the dense subspace ${\cal S}(\mb{R}^n)$ is closed under this action. For $\Om\in\mathscr{H}_n$ the corresponding Gaussian function $f_\Om(x)=e^{\pi i x\Om x^T}$ is in ${\cal S}(\mb{R}^n).$  The theta function $\t: \mathscr{H}_n\ra\mb{C}$ given by
\[\t(\Om)=\sum_{m\in\mb{Z}^n}f_\Om(m)=\sum_{m\in\mb{Z}^n}e^{\pi im\Om m^T}\] is automorphic for some arithmetic subgroup of $Sp_{2n}(\mb{R}).$ In \cite{Zhu} Y. Zhu has generalized
this classical theory of Weil representations and theta functions to the loop group $Sp_{2n}(\mb{R}((t))).$
His method also works in our situation. We do not attempt to build the theory of Weil representations for
our loop group in full generality. Instead we shall focus on theta functions. From now on we use row vectors instead of column vectors, and matrix group acts from the right.

Define $W=V\op V^\ast=\mb{R}[t,t^{-1}]^n\op \mb{R}[t,t^{-1}]^n$, where the canonical pairing
$(v,v^\ast)$ equals the constant term of $v(t)v^\ast(t)^T.$ The space $W$ has a symplectic form given by
$\lg v_1+v_1^\ast,v_2+v_2^\ast\rg=(v_1,v_2^\ast)-(v_2,v_1^\ast).$ We define $Sp_{2n}(\mb{R}[t,t^{-1}])$ to be the isometry group of $(W,\lg,\rg)$.

Denote by $V_\mb{C}$ and $V^\ast_\mb{C}$ the complexifications of $V$ and $V^\ast$; let $\mathscr{H}$ be the set of $\Om=X+iY\in \rm{Hom}_\mb{C}(V_\mb{C},V_\mb{C}^\ast)$ such that $X, Y\in \rm{Hom}_\mb{R}(V,V^\ast)$ are are self-dual and $Y>0.$ Note that $\rm{Im}(\Om)=Y$ induces an inner product $(w,v)_{iY}$ on $V$ via $(w,v)_{iY}:=(w,vY).$ Let $iI\in \mathscr{H}$ such that $I$ induces the inner product $(,)_1$ defined by (\ref{in}). In concrete terms, the operator $I: V\ra V^\ast$ maps $t^ie_j$ to $q^{-2i}t^{-i}e_j$, where $i\in\mb{Z}$ and $\{e_1,\ld,e_n\}$ is the standard basis of $\mb{R}^n.$  For any
\[
g=\bbm \a & \b \\ \g & \de\ebm\in Sp_{2n}(\mb{R}[t,t^{-1}]),
\]
$g$ acts on $\mathscr{H}$ by $g\c\Om=(\a\Om+\b)(\g\Om+\de)^{-1}.$ Consider the diagonal embedding
\[
d: GL_n(\mb{R}[t,t^{-1}])\hr
Sp_{2n}(\mb{R}[t,t^{-1}]),\q  g\mt d(g)=\bbm g^\ast & 0\\ 0 & g\ebm,
 \]
 where $g^\ast=(g^{-1})^T.$ Then $d(g)\c iI$ induces the inner product $(,)_g.$ This suggests that we may generalize the theta function to the loop symplectic group.

\begin{rmk} Here our diagonal embedding differs from the usual one because in previous sections we let $GL_n(V)$ act on $V$ from the left. For the setting of Weil representations by convention the symplectic group acts from the right.
\end{rmk}

We first generalize the notion of covolume theory. Define a sesquilinear form $(,)_\Om: V_\mb{C}\times V_\mb{C}\ra\mb{C}$ by $(v,w)_\Om=-i(\bar{v},w\Om).$ It is easy to check that $\rm{Re}(v,v)_\Om>0$ for any $v\neq 0.$ If $L\ss L'$ are two lattices of $R^n=\mb{Z}[t,t^{-1}]^n$,
the orthogonal complement $L_\mb{C}^\perp$ of $L_\mb{C}$ in $L'_\mb{C}$ with respect to the form $(,)_\Om$ is well-defined. If $\{v_i\}$ is a complementary basis  (i.e. $\{v_i\}$ span a complement of $L$ in $L'$), then we define the \textit{covolume} of $L$ in $L'$ to be the complex number
\be\label{cc}
c(\Om, L, L')=\sqrt{\det(v_i^\perp, v_j^\perp)}.
\ee
Here for the square root we take the positive branch. Similarly as in section 2, a \textit{covolume theory} for $\Om$ is a rule $c$ assigning a nonzero complex number $c(L)$ to each lattice $L$ of $R^n$ such that $c(L')=c(L)c(\Om, L, L')$ whenever $L\ss L'$ are two lattices.

If $\Om= g\c iI$ for some $g\in Sp_{2n}(\mb{R}[t,t^{-1}])$ and $c$ is a covolume theory for $\Om$, then we define
\be\label{tc}
\vt(\Om,c)=\lim_L c(L)\sum_{v\in L}e^{-\pi(v,v)_\Om}.
\ee

\begin{prop}\label{conc}
The limit defining $\vt(\Om,c)$ is convergent.
\end{prop}

The proof of Proposition \ref{conc} relies on the Iwasawa decomposition for loop symplectic groups, which we shall recall now. Let
\be \label{ksp}
K_{sp}=\{g\in Sp_{2n}(\mb{R}[t,t^{-1}])|g(q^{-1}t)g(q^{-1}t^{-1})^T=1\}.
\ee
It can be shown that $K_{sp}$ is the stabilizer of $iI$ in $Sp_{2n}(\mb{R}[t,t^{-1}]).$ Let $B_0$ denote the Borel subgroup of $Sp_{2n}(\mb{R})$ given by
\be
B_0=\l\{\l.\bbm A^\ast & C\\ 0 & A\ebm\in Sp_{2n}(\mb{R})\r| A\rm{ is upper triangular}\r\},
\ee
and define
\be
B_{sp}=\{g\in Sp_{2n}(\mb{R}[t^{-1}])| g(\i)\in B_0\}.
\ee
Then we have the Iwasawa decomposition
\be
Sp_{2n}(\mb{R}[t,t^{-1}])=B_{sp}K_{sp}.
\ee

\

\textit{Proof} of Proposition \ref{conc}: It suffices to consider $\Om=g\c iI$ with $g\in B_{sp}$. We use the fact that $B_{sp}$ is generated by elements of the form
\be\label{gen}
\bbm \a & 0\\ \g & \de\ebm\rm{\q or\q } \bbm \a & \b \\ 0 & \de\ebm.
\ee
Note that in this case $\de=\a^\ast\in GL_n(\mb{R}[t^{-1}])$, $\b, \g\in M_{n\times n}(\mb{R}[t^{-1}]).$ Moreover $\de\in B.$ We may write $g=g_l\cs g_1$ such that each $g_j$ is of the form (\ref{gen}).
Hence we have a chain of transformations
\[
iI=\Om_0\stackrel{g_1}{\lra} \Om_1\stackrel{g_2}{\lra}  \Om_2\stackrel{g_3}{\lra} \cs \stackrel{g_l}{\lra}\Om_l=\Om. \]
For each $j$ one has either $\Om_j=\a_j\Om_{j-1}\a^T_j+\b_j\a_j^T$ or $\Om_j=(\g_j\de_j^T+\de_j\Om^{-1}_{j-1}\de_j^T)^{-1}$. One may keep track of these finitely many steps and apply the same method but with some variant of the arguments used in Section 3.\hfill$\Box$

\section{Interpretations of theta functions}

We start from general theta functions with parameters. Keep the notations in the last section. Let $\Om=g\c iI\in\mathscr{H}$, where $g\in Sp_{2n}(\mb{R}[t,t^{-1}])$. Let $c$ be a covolume theory for $\Om.$ For $a\in V$, $b\in V^\ast$, $z\in V^\ast+V\Om= V^\ast_\mb{C}=\mb{C}[t,t^{-1}]^n$, we define the theta function
\be\label{gt}
\vt \bbm a\\ b\ebm(z,\Om,c)=\lim_L c(L)\sum_{v\in L}\exp\l(-\pi(v+a,v+a)_\Om+2\pi i(v+a,z+b)\r).
\ee
We have the quasi-periodicity: for $v\in\mb{Z}[t,t^{-1}]^n$,
\bq
&& \vt\bbm a\\ b\ebm(z+v,\Om,c)=e^{2\pi i(a,v)}\vt\bbm a\\ b\ebm(z,\Om,c),\label{per1}\\
&& \vt\bbm a\\ b\ebm(z+v\Om,\Om,c)=e^{\pi(v,v)_\Om-2\pi i(v,z+b)}\vt\bbm a\\ b\ebm(z,\Om,c).\label{per2}
\eq
Let $L_\Om=\mb{Z}[t,t^{-1}]^n+\mb{Z}[t,t^{-1}]^n\Om\ss V^\ast_\mb{C}=\mb{C}[t,t^{-1}]^n.$
Let $l$ be a positive integer and fix a covolume theory $c_l$ for $l\Om$. Each $a\in \mb{Z}[t,t^{-1}]^n$ gives a theta function
\be\label{fa}
f^a_l(z):=\vt\bbm a/l \\ 0\ebm (lz,l\Om,c_l).
\ee
 For a lattice $L$ of $\mb{Z}[t,t^{-1}]^n$, introduce the partial theta function
\be\label{parth}
f^a_{l,L}(z)=\vt_L\bbm a/l \\\ 0\ebm (lz, l\Om, c_l).
\ee
From quasi-periodicity (\ref{per1}), (\ref{per2}) it follows that for $v\in \mb{Z}[t,t^{-1}]^n$,
\bq
&& f^a_l(z+v)=f^a_l(z),\label{per3}\\
&& f^a_l(z+v\Om)=e^{\pi l(v,v)_\Om-2\pi i l(v,z)}f^a_l(z).\label{per4}
\eq
Following \cite{Zhu}, let $S_l$ be the set of $a\in \mb{Z}[t,t^{-1}]^n$ such that all the coefficients in $a$ lie in $\{0,1,\ld,l-1\}.$ Let
\[
\mathscr{B}_l=\{z\in V_\mb{C}^\ast|f^a_l(z)=0\rm{ for all }a\in S_l\}/L_\Om,
\]
which is called the set of ``base points" in the complex torus $V_\mb{C}^\ast/L_\Om.$ Let $\mb{P}^{S_l}=\mb{P}(\mb{C}^{S_l})$ be the infinite projective space, then there is a canonical holomorphic map
\[
\phi_l: V_\mb{C}^\ast/L_\Om-\mathscr{B}_l\ra \mb{P}^{S_l}
\]
given by $\phi_l(z)=[f^a_l(z)]_{a\in S_l}.$ The well-definedness of $\phi_l$ follows from quasi-periodicity (\ref{per3}), (\ref{per4}). It is known that $\mathscr{B}_l=\emptyset$ for $l\geq 2$ and $\phi_l$ is an embedding for $l\geq 3.$

We introduce a line bundle $\mathscr{L}$ over $V_\mb{C}^\ast/L_\Om.$ Define an action of $L_\Om$ on the trivial bundle $V_\mb{C}^\ast\times\mb{C}\ra V_\mb{C}^\ast$ by
 \[
 (v_1+v_2\Om)\c(z,\la)=(z+v_1+v_2\Om, e^{\pi(v_2,v_2)\Om-2\pi i (v_2,z)}\la).
 \]
Let $\mathscr{L}$ be the quotient line bundle. If $f$ is a function on $V^\ast_\mb{C}$ satisfying (\ref{per3}) (\ref{per4}), and $f$ is holomorphic, i.e. the restriction of $f$ to any finite dimensional complex subspace of $V_\mb{C}^\ast$ is  holomorphic, then we say $f$ is a global section of $\mathscr{L}^l$. Let
$\G(\mathscr{L}^l)$ be the set of global sections of $\mathscr{L}^l.$

We claim that $f^a_l$, $a\in S_l$ are  global sections of $\mathscr{L}^l$. It suffices to prove the holomorphic property. This follows from certain uniform convergence (in the variable $z$) results, and one may argue similarly as before. We omit the details. If we fix a covolume theory $c_l'$ for $l^{-1}\Om$ and for $b\in S_l$ define
\[
g^b_l(z)=\vt\bbm 0\\ b/l\ebm (z,l^{-1}\Om,c_l'),
\]
then $g^b_l$, $b\in S_l$ are also global sections of $\mathscr{L}^l.$ In the special case $l=k^2$, there is another family of global sections
\[
h^{a,b}_l(z)=\vt\bbm a/k\\ b/k\ebm (kz,\Om,c),\q a, b\in S_k.
\]

Let us introduce the loop Heisenberg group $Heis(n,l)$, which is the set $U_l\times \mb{Z}_l[t,t^{-1}]^n\times\mb{Z}_l[t,t^{-1}]^n$ subject to the following relations:
\bq
&& (1,x,0)(1,x',0)=(1,x+x',0),\\
&& (1,0,y)(1,0,y')=(1,0,y+y'),\\
&& (\zeta,0,0)(\zeta',0,0)=(\zeta\zeta',0,0), \rm{ and }U_l\times\{0\}\times\{0\}\rm{ is the center},\\
&& (1,0,y)(1,x,0)=\zeta_l((x,y))(1,x,0)(1,0,y),
\eq
where $(x,y)\in\mb{Z}_l$ is the constant term of $x(t)y(t)^T$, and $\zeta_l(a):=(\zeta_l,a)$, $a\in \mb{Z}_l$, is the canonical pairing.

$\G(\mathscr{L}^l)$ can be realized as a representation of $Heis(n,l)$ as follows. Note that there is a natural bijection between $S_l$ and $\mb{Z}_l[t,t^{-1}]^n.$ For $x\in \mb{Z}_l[t,t^{-1}]^n$ write $\wtl{x}$ for the corresponding element in $S_l.$ Conversely for $a\in S_l$ write $\bar{a}=a\mod l\in\mb{Z}_l[t,t^{-1}]^n.$ We define an action of $Heis(n,l)$ on $f(z)\in \G(\mathscr{L}^l)$ by
\bq
&& (1,x,0)f(z)=f(z+\wtl{x}/l),\\
&& (1,0,y)f(z)=\exp(-\pi(\wtl{y},\wtl{y})_\Om/l+2\pi i (\wtl{y},z))f(z+\wtl{y}\Om/l),\\
&& (\zeta,0,0)f(z)=\zeta f(z).
\eq
One can check that this defines a group action of $Heis(n,l)$ on $\G(\mathscr{L}^l).$ Acting on above three families of sections we have the formulas
\[
\ba{ll}
(1,x,0)f^a_l=\zeta_l((\bar{a},x))f^a_l, & (1,0,y)f^a_l=f_l^{a+\wtl{y}},\\
(1,x,0)g^b_l=g_l^{b+\wtl{x}},  & (1,0,y)g^b_l=\zeta_l(-(y,\bar{b}))g^b_l,\\
(1,x,0)h_l^{a,b}=h_l^{a,b+k\wtl{x}}, & (1,0,y)h_l^{a,b}=\zeta_k(-(y,\bar{b}))h_l^{a+k\wtl{y},b}.
\ea
\]

Consider the graded algebra $\bop\limits^\i_{l=0}\G(\mathscr{L}^l)$, where $\G(\mathscr{L}^0):=\mb{C}$. We shall give the explicit formula for $f_{l_1}^{a_1}f_{l_2}^{a_2}$, where $a_i\in S_{l_i}$, $i=1,2$, under the multiplication $\G(\mathscr{L}^{l_1})\times \G(\mathscr{L}^{l_2})\ra \G(\mathscr{L}^{l_3})$, where $l_3=l_1+l_2.$ Intuitively we may expand $f_{l_1}^{a_1}f_{l_2}^{a_2}$ as an infinite linear combination of $f_{l_3}^a$, $a\in S_{l_3}$. However it turns out to be not the case. The formalism of the result is by first averaging partial theta functions and then taking the limit.

We need to introduce some notations for our formula. Let $(l_1,l_2)$ be the greatest common divisor of $l_1$ and $l_2$, and let $l_i'=l_i/(l_1,l_2)$, $i=1,2,3.$ Choose $j_1, j_2\in\mb{Z}$ such that $j_1l_1+j_2l_2=(l_1,l_2).$ For any $\eta\in\mb{Z}[t,t^{-1}]^n$, define
 \bq
 &&a_\eta=a_1+a_2+(l_1,l_2)\eta,\\
 && \wtl{a}_\eta=l_2'a_1-l_1'a_2+l_1l_2'(j_1-j_2)\eta.
 \eq
Note that $(\wtl{a}_\eta\mod l_1'l_2'l_3)$ does not depend on the choice of $(j_1, j_2)$.
Introduce a subset of $\mb{Z}[t,t^{-1}]^n$,
\be
S^{a_1,a_2}_{l_1,l_2}=\{\eta\in\mb{Z}[t,t^{-1}]^n|a_\eta\in S_{l_3}\}.
\ee
It is clear that $S^{a_1,a_2}_{l_1,l_2}$ is a set of representatives for $\mb{Z}_{l_3'}[t,t^{-1}]^n$, i.e.
\be
\mb{Z}[t,t^{-1}]^n=\bigsqcup_{\eta\in S^{a_1,a_2}_{l_1,l_2}}\eta+l_3'\mb{Z}[t,t^{-1}]^n.
\ee
Similarly, for $d\in\mb{Z}$ such that $a_1, a_2\in L^d_\mb{Z}$ one has
\be\label{dis}
L^d_\mb{Z}=\bigsqcup_{\eta\in S^{a_1,a_2}_{l_1,l_2}\cap L^d_\mb{Z}}\eta+l_3'L^d_\mb{Z}.
\ee

We also need some notions from \cite{Kap}. Let us denote by $\mathscr{G}(V)$ the set of lattices of $V$ and call it the \textit{semiinfinite Grassmannian} of $V$. A \textit{dimension theory} on $V$ is a map $\mathscr{D}: \mathscr{G}(V)\ra\mb{Z}$ such that whenever $L, L'\in \mathscr{G}(V)$ and $L\ss L'$, we have
$\mathscr{D}(L')=\mathscr{D}(L)+\dim (L'/L).$ All dimension theories on $V$ form a $\mb{Z}$-torsor.

Recall that we have fixed a covolume theory $c_j$ for $j\Om$, $\forall j>0.$ We claim that there exist a dimension theory $\mathscr{D}$ on $V$ and a constant $\la\in \mb{C}^\times$ such that
\be\label{dim}
c_{l_1}c_{l_2}= \la\c( l'_3)^{-\mathscr{D}}\c c_{l_1'l_2'l}c_{l_3}.
\ee
Note that we view both sides of (\ref{dim}) as functions $\mathscr{G}(V)\ra\mb{C}^\times.$ The claim follows from the simple identity
\[
l_1l_2=l'_1l_2'l^2_3 (l'_3)^{-2}.
\]
The choice of $(\mathscr{D}, \la)$ is not unique, and we only need to fix such a pair. Now we are ready to give the
asymptotic multiplication formula.

\begin{prop}\label{mul}
With above notations, and a choice of the pair $(\mathscr{D}, \la)$ satisfying $($\ref{dim}$)$, we have
\[
f_{l_1}^{a_1}(z)f_{l_2}^{a_2}(z)=\la\c\lim_{d\ra\i} (l_3')^{-\mathscr{D}(L^d_\mb{R})}\sum_{\eta\in S^{a_1,a_2}_{l_1,l_2}\cap L^d_\mb{Z}}f^{\wtl{a}_\eta}
_{ l_1'l_2'l_3, L^d_\mb{Z}}(0)\c f^{a_\eta}_{l_3,L^d_\mb{Z}}(z).
\]
\end{prop}

\begin{cor}\label{tre}
If $l_1=l_2=l$, then
\[
f^{a_1}_l(z)f^{a_2}_l(z)=\la\c\lim_{d\ra\i} 2^{-\mathscr{D}(L^d_\mb{R})}\sum_{\eta\in S_2}f^{a_1-a_2+l\eta}_{2l, L^d_\mb{Z}}(0)f^{a_1+a_2+l\eta}_{2l,L^d_\mb{Z}}(z).
\]
\end{cor}

\begin{rmk} Corollary \ref{tre} is an infinite analog of Corollary 6.8 in \cite{Mum}, which plays a major role in Riemann's theta relation $($see \cite{Mum} \S7$)$.
\end{rmk}

\textit{Proof} of Proposition \ref{mul}: Assume that $a_1, a_2\in L^d_\mb{Z}$ such that (\ref{dis}) is satisfied. We shall compute $f^{a_1}_{l_1,L^d_\mb{Z}}(z)f^{a_2}_{l_2, L^d_\mb{Z}}(z)$ explicitly.
Consider the following expression for $v_1, v_2$ run over $L^d_\mb{Z}\times L^d_\mb{Z}$:
\bq\label{sum}
&&l_1(v_1+\f{a_1}{l_1},v_1+\f{a_1}{l_1})_{\Om}-2i(v_1+\f{a_1}{l_1}, l_1z)\\
&+&l_2(v_2+\f{a_2}{l_2},v_2+\f{a_2}{l_2})_{\Om}-2i(v_2+\f{a_2}{l_2}, l_2z).\n
\eq
Make the change of variables
\[
v_1=j_1(\eta+l_3'v)+l_2'u,\q v_2=j_2(\eta+l_3'v)-l_1'u,
\]
such that $\eta$ runs over $S^{a_1,a_2}_{l_1,l_2}\cap L^d_\mb{Z}$, and $u$, $v$ run over $L^d_\mb{Z}.$ Then from (\ref{dis}) it follows that
\[
L^d_\mb{Z}\times L^d_\mb{Z}\ra (S^{a_1,a_2}_{l_1,l_2}\cap L^d_\mb{Z})\times L^d_\mb{Z}\times L^d_\mb{Z}, \q (v_1, v_2)\mt (\eta, u, v)
\]
is a bijection. In terms of the new variables, after some heavy manipulations one can check that (\ref{sum}) equals
\bq
&&l_1'l_2'l_3(u+\f{\wtl{a}_\eta}{l_1'l_2'l_3}, u+\f{\wtl{a}_\eta}{l_1'l_2'l_3})_\Om\n\\
&+& l_3(v+\f{a_\eta}{l_3},v+\f{a_\eta}{l_3})_\Om-2i(v+\f{a_\eta}{l_3}, l_3z).\n
\eq
Taking into account of (\ref{dim}), this finishes the proof.\hfill$\Box$

\begin{rmk}
 When $q\in\mb{Q}$, the endomorphism ring of $V^\ast_\mb{C}/L_\Om$ is much larger than $\mb{Z},$ namely, it contains $\mb{Z}[Nt, Nt^{-1}]$ for some sufficiently large integer $N$. This is an interesting analogue of abelian varieties with large endomorphism rings.
\end{rmk}

\section*{Acknowledgement}

The author is grateful to Prof. Y. Zhu for very stimulating and helpful discussions during this work.

\bibliographystyle{amsplain}

\end{document}